\documentclass[12pt]{article}
\usepackage{amssymb,amsmath,graphicx,verbatim,eufrak,amsthm}
\usepackage{hyperref}
\newtheorem{thm}{Theorem}[section]

\newtheorem*{conj*}{Conjecture}

\newtheorem{cor}[thm]{Corollary}
\newtheorem{clm}[thm]{Claim}
\newtheorem{lem}[thm]{Lemma}
\newtheorem{prop}[thm]{Proposition}
\theoremstyle{definition}

\theoremstyle{remark}
\newtheorem{rem}[thm]{Remark}

\newcommand{\Claim}{\textbf{Claim.}\hspace{10pt}}

\renewcommand{\P}{\mathcal{P}}

\newcommand{\Rng}{\mathbf{Rng}}

\newcommand{\abs}[1]{\left\vert#1\right\vert}

\newcommand{\set}[1]{\left\{#1\right\}}
\newcommand{\setb}[2]{ \left\{#1 \ \big| \ #2 \right\} }

\newcommand{\Set}[2]{ \left\{#1 \ \big| \ #2 \right\} }
\newcommand{\br}[1]{\left[#1\right]}

\newcommand{\sr}[1]{\left(#1\right)}

\newcommand{\Bracket}[1]{\left[#1\right]}
\newcommand{\Soger}[1]{\left(#1\right)}

\newcommand{\Integer}{\mathbb{Z}}

\newcommand{\eps}{\varepsilon}

\newcommand{\E}{\mathbf{E}}

\def\squareforqed{\hbox{\rlap{$\sqcap$}$\sqcup$}}
\def\qed{\ifmmode\squareforqed\else{\unskip\nobreak\hfil
\penalty50\hskip1em\null\nobreak\hfil\squareforqed
\parfillskip=0pt\finalhyphendemerits=0\endgraf}\fi}

\newcommand{\bgive}[2]{\Bracket{#1 \ \big| \ #2}}

\renewcommand{\th}{^{\mathrm{th}}}

\renewcommand{\H}{\mathcal{H}}

\newcommand{\Hom}{\mathrm{Hom}}
\newcommand{\Lip}{\mathrm{Lip}}

\newcommand{\Z}{\Integer}
\newcommand{\N}{\mathbb{N}}

\newcommand{\dwn}{\downarrow}

\renewcommand{\O}{\mathcal{O}}

\newcommand{\eqdef}{\stackrel{\tiny \textrm{def} }{=}}
\begin{document}

\author{Itai Benjamini \\
Department of Mathematics \\
Weizmann Institute \\
Rehovot 76100, Israel \\
\texttt{itai.benjamini@weizmann.ac.il} %
\and Ariel Yadin \\
Department of Mathematics \\
Weizmann Institute \\
Rehovot 76100, Israel \\
\texttt{ariel.yadin@weizmann.ac.il} %
\and Amir Yehudayoff \\
Department of Computer Science \\
Weizmann Institute \\
Rehovot 76100, Israel \\
\texttt{amir.yehudayoff@weizmann.ac.il} }

\title{Random Graph-Homomorphisms and Logarithmic Degree}

\date{ }

\maketitle


\begin{abstract}
A graph homomorphism between two graphs is a map from the vertex
set of one graph to the vertex set of the other graph, that maps
edges to edges. In this note we study the range of a uniformly
chosen homomorphism from a graph $G$ to the infinite line $\Z$. It
is shown that if the maximal degree of $G$ is `sub-logarithmic',
then the range of such a homomorphism is super-constant.

Furthermore, some examples are provided, suggesting that perhaps
for graphs with super-logarithmic degree, the range of a typical
homomorphism is bounded. In particular, a sharp transition is
shown for a specific family of graphs $C_{n,k}$ (which is the
tensor product of the $n$-cycle and a complete graph, with
self-loops, of size $k$). That is, given any function $\psi(n)$
tending to infinity, the range of a typical homomorphism of
$C_{n,k}$ is super-constant for $k= 2\log(n) - \psi(n)$, and is
$3$ for $k= 2\log(n) + \psi(n)$.
\end{abstract}




Submitted to EJP on November 29 2006, final version accepted June
11 2007.


\section{Introduction}


A graph homomorphism from a graph $G$ to a graph $H$ is a map from
the vertex set of $G$ to the vertex set of $H$, that maps edges to
edges. By a \emph{homomorphism of $G$} we mean a graph
homomorphism from $G$ to the infinite line $\Z$. Thus, a
homomorphism of $G$ maps adjacent vertices to adjacent integers.
We note that the uniform measure on the set of all homomorphisms
of $G$, that send some fixed vertex to $0$, generalizes the
concept of random walks on $\Z$. Indeed, a random homomorphism of
the $k$-line is a random walk of length $k$ on $\Z$. So, random
homomorphisms of a graph $G$, are also referred to as $G$-indexed
random walks. Tree-indexed random walks were studied by Benjamini
and Peres in \cite{BP}.  For results concerning random
homomorphisms of general graphs see \cite{BHM, LNR}. \cite{BS}
deals with connections between random homomorphisms and the
Gaussian random field.  For other related 2-dimensional height
models in physics see \cite{Kenyon, Sheffield}.

A key quantity for our first result, Theorem \ref{lower bound
theorem}, is $V(r)$, the maximal size of a ball of radius $r$ in
$G$. Theorem \ref{lower bound theorem} states that for every $r$
such that $V(r)$ is at most $\frac{1}{2} \log(\abs{G})$, the range
of a random homomorphism is greater than $r$, with high
probability.  If $d$ is the maximal degree in $G$, then $V(r)$ is
at most $(d+1)^{r}$. Thus, Theorem \ref{lower bound theorem}
implies that for graphs of `small enough' degree, the range of a
homomorphism is typically `large' (see Corollary~\ref{lower bound
corrolary}). We stress that this is only a sufficient condition
for large range, and not a necessary one. For example, consider
the $\log(n)$-regular tree of size $n$.  Already for $r=1$, a ball
of radius $r$ has at least $\log(n)$ vertices, so the assertion of
Theorem \ref{lower bound theorem} is trivial. However, the range
of a typical homomorphism of this tree is of size at least
$\Omega(\log(n)/\log\log(n))$.

The next natural question is: How tight is this lower bound?  That
is, are there examples of graphs of logarithmic degree that have
bounded range (as the size of the graph grows to infinity)?  This
can be decided via a result of Kahn \cite{Kahn}.  Kahn's results
states that there exists a constant $b \in \N$, such that the
range of a random homomorphism of $Q_d$, the discrete cube of
dimension $d$, is at most $b$, with probability tending to $1$ as
$d$ tends to infinity (note that the size of $Q_d$ is $2^d$ and
the degree of $Q_d$ is $d$). Galvin \cite{Galvin} later calculated
$b=5$.

Kahn's result raises a new question: What happens to the range of
a random homomorphism, if the degree is logarithmic, but the
diameter is large? (The discrete cube has logarithmic degree, but
also has logarithmic diameter.) To answer this question, we study
the graph $C_{n,k}$ in Section \ref{sec: the cycle} (the graph
$C_{n,k}$ is the tensor product of the $n$-cycle and the complete
graph of size $k$ with self-loops). We show a sharp transition in
$k$, of the range of a random homomorphism of $C_{n,k}$. Namely,
for any monotone function $\psi(n)$ tending to infinity, if $k =
2\log(n) - \psi(n)$ the range is $2^{ \Omega (\psi(n))}$, with
high probability, and if $k = 2\log(n) + \psi(n)$ the range is
$3$, with high probability. In particular, $C_{n,3\log n}$ is a
graph of almost linear diameter and logarithmic degree such that
the range of a random homomorphism of $C_{n,3 \log n}$ is $3$,
with high probability.

The rest of this paper is organized as follows:  We first
introduce some notation.  Section \ref{sec: lower bound} contains
our lower bound.  Section \ref{sec: the cycle} proves the upper
and lower bounds on the range of random homomorphisms of the graph
$C_{n,k}$.  Section \ref{sec: further research} lists some further
possible research directions concerning random homomorphisms of
graphs.

{\bf Acknowledgement.} \hspace{10pt} We would like to thank Ori
Gurel-Gurevich for useful discussions.

\subsection{Notation and Definitions}

Logarithms are always of base $2$. For an integer $k \in \N$,
denote $[k] = \set{1,\ldots,k}$. For two integers $x,y \in \Z$,
denote by $[x,y]$ the set of integers at least $x$ and at most
$y$. For $n \in \N$, denote by $\Z_n$ the additive group whose
elements are $[0,n-1]$, and addition is modulo $n$.

\subsubsection{Graphs}

All graphs considered are simple and connected. Let $G$ be a
graph. For simplicity of notation, we use $G$ to denote the vertex
set of the graph $G$. In particular, we write $v \in G$, if $v$ is
a vertex of the graph $G$. For two vertices $v,u \in G$, we write
$v \sim u \in G$, if $\set{u,v}$ is an edge in the
graph $G$. When the graph is clear, we use $v \sim u$. 
The \emph{size} of the graph $G$, denoted $\abs{G}$, is the number
of vertices in $G$. The \emph{diameter} of $G$ is the maximal
distance between any two vertices in $G$. For a vertex $v \in G$
and an integer $r \in \N$, a \emph{ball of radius $r$ centered at
$v$} is the subgraph of $G$ induced by the set of all vertices at
distance at most $r$ from $v$.

\subsubsection{Homomorphisms}

For two graphs $G$ and $H$, a \emph{graph homomorphism from $G$ to
$H$} is a mapping $f: G \to H$ that preserves edges; that is,
every two vertices $v \sim u \in G$ admit $f(v) \sim f(u) \in H$.
For two vertices $v_0 \in G$ and $x_0 \in H$, we denote by
$\Hom_{v_0}^{x_0}(G,H)$ the set of all homomorphisms $f$ from $G$
to $H$ such that $f(v_0) = x_0$. A homomorphism from $G$ to $H$ is
also called a \emph{$H$-coloring of $G$}.

We denote by $\Z$ both the set of integers, and the graph whose
vertex set is the integers and edge set is $\setb{\set{z,z+1}}{z
\in \Z}$. We mostly consider $\Hom_{v_0}(G,\Integer) \eqdef
\Hom_{v_0}^{0}(G,\Integer)$. Note that
$$\Hom_{v_0}(G, \Z) =
\setb{ f:G \to \Z}{ \forall \ u \sim v \in G \ \abs{f(u) - f(v)} =
1 \textrm { and } f(v_0) = 0 } . $$ %
For a mapping $f:G \to \Z$, define
$$ f(G) = \setb{ f(v)}{v \in G} \quad \textrm{ and } \quad R(f) =
\abs{f(G)} . $$ %
We call both $f(G)$ and $R(f)$ the \emph{range} of $f$. We use the
notation $\in_R$ to denote an element chosen uniformly at random.
E.g., $f \in_R \Hom_{v_0}(G,\Integer)$ is a random homomorphism
from $G$ to $\Integer$ such that $f(v_0) = 0$, chosen uniformly at
random. (When $G$ is finite and connected, the set
$\Hom_{v_0}(G,\Z)$ is finite, and $f \in_R \Hom_{v_0}(G,\Integer)$
is well defined.) For example, consider the case where $G$ is the
interval of length $n$; that is
$$ V(G) = [0,n] \quad \textrm{ and } \quad E(G) = \setb{ \set{i,i+1}
}{ 0 \leq i \leq n-1 } . $$ %
Then, $\Hom_0(G,\Z)$ is the set of all paths in $\Z$ starting from
$0$, of length $n$. Therefore, $f \in_R \Hom_0(G,\Z)$ is a
$n$-step random walk on $\Z$, starting at $0$.  Thus, for a
general (connected and finite) graph $G$, a random homomorphism $f
\in_R \Hom_{v_0}(G,\Z)$, is also called a \emph{$G$-indexed random
walk}.

For a graph $G$, we say that a homomorphism $f$ from $G$ to itself
is an \emph{automorphism}, if $f$ is invertible, and $f^{-1}$ is a
homomorphism as well. We say that a graph $G$ is \emph{vertex
transitive}, if all the vertices of $G$ ``look'' the same; that
is, for any two vertices $v,u \in G$, there exists an automorphism
$f$ of $G$ such that $f(v) = u$. We say that a graph $G$ is
\emph{edge transitive}, if all the edges of $G$ ``look'' the same;
that is, for any two edges $\set{v_1,v_2}$ and $\set{u_1,u_2}$ in
$G$, there exists an automorphism $f$ of $G$ such that
$\set{f(v_1),f(v_2)} = \set{u_1,u_2}$.

\section{Lower Bounds for Graphs with Small Degree} \label{sec:
lower bound}

In this section we show that for graphs of `small enough' degree,
the range of a homomorphism is typically `large'. In fact, we
prove something slightly stronger:

\begin{thm} \label{lower bound theorem}
Let $\set{G_n}$ be a family of graphs such that $\lim_{n \to
\infty} \abs{G_n} = \infty$. For $r \in \N$, define $V_n(r)$ to be
the maximal size of a ball of radius $r$ in $G_n$. Let $v_n \in
G_n$ and let $f_n \in_R \Hom_{v_n}(G_n,\Z)$ be a random
homomorphism. Let $r = r(n) \in \N$. Assume that there exists a
constant $c < 1$ such that every large enough $n \in \N$ admits
$V_n(r) \leq c \log \abs{G_n}$. Then
$$ \Pr \Bracket{ R(f_n) \leq r} = o(1) . $$
\end{thm}

We defer the proof of Theorem~\ref{lower bound theorem} to
Section~\ref{sec: proof of lower bound gen}. First we discuss the
tightness of Theorem~\ref{lower bound theorem}. In
Section~\ref{sec: the cycle} we consider the family of graphs
$\set{C_{n,k}}$, where $n \in \N$ is even, and $k = k(n) \in \N$.
For $n \in \N$, the size of $C_{n,k}$ is $kn$, and the size of a
ball of radius $3$ in $C_{n,k}$ is at most $7k$; that is, $V_n(3)
\leq 7k$. In Theorem~\ref{thm: upper bound for cnk} we show an
upper bound on the range of a random homomorphism of $C_{n,k}$,
for logarithmic $k$. More specifically, we show that for $k =
2\log n + \log \log \log n$,
\begin{eqnarray} \Pr \Bracket{R(f_n) > 3} = o(1), \label{eqn: nn1} \end{eqnarray} where $f_n \in_R
\Hom_{(0,1)}(C_{n,k},\Z)$ is a random homomorphism. Thus,
Theorem~\ref{lower bound theorem} is wrong if instead of $c<1$ we
use $c \leq 14$. Indeed, assume towards a contradiction that
Theorem~\ref{lower bound theorem} holds for every $c \leq 14$.
Since $V_n(3) \leq 7k < 14 \log (kn)$, by the assumption we have
$$ \Pr \Bracket{ R(f_n) \leq 3 } = o(1),$$
where $f_n \in_R \Hom_{(0,1)}(C_{n,k},\Z)$ is a random
homomorphism. This is a contradiction to (\ref{eqn: nn1}). We note
that since $C_{n,k}$ is vertex transitive (and edge transitive),
Theorem~\ref{lower bound theorem} is tight in the above sense for
vertex transitive graphs (and for edge transitive graphs).

\subsection{Lower Bounds for Graphs with Small Degree}

The following corollary of Theorem \ref{lower bound theorem} shows
that the range of a random homomorphism from a graph of ``small''
degree to $\Z$ is ``large''. For example, consider any family of
graphs $\set{G_n}$, such that the degree of $G_n$ is $\log \log
\abs{G_n}$.  Then, the corollary states that the range of a random
homomorphism from $G_n$ to $\Z$ is super-constant (as $n$ tends to
infinity), with high probability.
\begin{cor} \label{lower bound corrolary}
Let $\set{G_n}$ be a family of graphs such that $\lim_{n \to
\infty} \abs{G_n} = \infty$. Let $n \in \N$, and let $d=d(n)$ be
the maximal degree of $G_n$. Let $v_n \in G_n$ and let $f_n \in_R
\Hom_{v_n}(G_n,\Z)$ be a random homomorphism. Then
$$ \Pr \Bracket{ R(f_n) \leq \frac{\log \log \abs{G_n} -1}{\log(d+1)} } = o(1) . $$
\end{cor}

\begin{proof}
For $r \in \N$, denote by $V_n(r)$ the maximal size of a ball of
radius $r$ in $G_n$. Since the maximal degree of $G_n$ is $d =
d(n)$, every $r \in \N$ admits $V_n(r) \leq (d+1)^r$. Denote $r =
r(n) = \frac{\log \log \abs{G_n} -1}{\log(d+1)}$. Since $(d+1)^r =
\frac{1}{2} \log \abs{G_n}$, the corollary follows, by Theorem
\ref{lower bound theorem} (with $c = 1/2$).
\end{proof}

\subsection{An Example - the Torus}

A specific example for the use Theorem~\ref{lower bound theorem}
is in the case of the torus. For an integer $n \in \N$, define the
$n \times n$ torus, denoted $T_n$, as follows: The vertex set is
$\Z_n \times \Z_n$, and the edge set is defined by the relations
$$\forall \ i,j \in \Z_n \ \ (i,j) \sim (i+1,j) \ \ \ (i,j) \sim (i,j+1)$$
(where addition is modulo $n$). Note that $T_n$ is both vertex
transitive and edge transitive. The following corollary shows that
the range of a random homomorphism of the $n \times n$ torus is at
least $\Omega( \log^{1/2} n)$, with high probability.

\begin{cor} Let $n \in \N$, and let $T_n$ be the $n \times n$ torus. Let $f_n \in_R
\Hom_{(0,0)}(T_n,\Z)$ be a random homomorphism. Then
$$\Pr[R(f_n) > 1/2 \log^{1/2} n] = 1-o(1).$$
\end{cor}

\begin{proof}
Note that the size of $T_n$ is $n^2$. For $r \in \N$, denote by
$V_n(r)$ the maximal size of a ball of radius $r$ in $T_n$. The
vertex set of the ball of radius $r$ centered at $(0,0)$ is
contained in
$$\setb{(i,j) \in T_n}{i,j \in \set{-r,\ldots,0,\ldots,r} \textrm{(modulo $n$)} } . $$
Thus, since $T_n$ is vertex transitive, $V_n(r) \leq (2r+1)^2$.
So, since every large enough $n$ admits $(2 (1/2 \log^{1/2} n)
+1)^2 \leq 2/3 \log (n^2)$, using Theorem~\ref{lower bound
theorem} (with $r = 1/2 \log^{1/2} n$ and $c = 2/3$), we have
$\Pr[R(f_n)
> r] = 1-o(1)$.
\end{proof}

\subsection{A Ball of Radius $r$ Has a Homomorphism with Range $r+1$}

Before we prove Theorem \ref{lower bound theorem}, we need the
observation of Lemma \ref{boundary}. 

%

For a graph $G$ and a vertex $v \in G$, we denote by $B_r(v)$ the
ball of radius $r$ centered at $v$. We say that \emph{$B_r(v)$ is
of exact radius  $r$}, if there exists $u \in B_r(v)$ such that
the distance between $v$ and $u$ in $G$ is at least $r$.

\begin{lem} \label{boundary} Let $G$ be a graph,
let $v$ be a vertex in $G$, and let $r \in \N$. For an integer $s
\in [0,r]$, define $B_s = B_s(v)$ to be the ball centered at $v$
with radius $s$. Set $B=B_r(v)$ and $\Gamma = \setb{u \in B}{u
\not \in B_{r-1}}$ (the boundary of $B$). Let $f$ be a
homomorphism from $B$ to $\Z$. Assume that $B$ is of exact radius
$r$. Then there exists a homomorphism $g$ from $B$ to $\Z$ such
that $g|_{\Gamma} = f|_{\Gamma}$, and $R(g) \geq r+1$.
\end{lem}

\begin{proof}
Since a translation of a homomorphism is a homomorphism with the
same range size, assume without loss of generality that $\min_{u
\in \Gamma} f(u) = 0$. We demonstrate an iterative process such
that at the $i\th$ step (for $i = 0,\ldots,r$) we have a
homomorphism $g_i$ that admits
\begin{enumerate}
    \item $g_i|_{\Gamma} = f|_{\Gamma}$.
    \item The minimal value of $g_i$ on the ball $B_{r-i}$ is
    $i$.
\end{enumerate}
Thus, setting $g = g_r$, we have: By property 1, we have
$g|_{\Gamma} = f|_{\Gamma}$. By property 2, we have $\max_{u \in
B} g(u) \geq r$. Thus, since $\min_{u \in \Gamma} g(u) = \min_{u
\in \Gamma} f(u) = 0$, we have $R(g) \geq \abs{[0,r]} = r+1$
(which completes the proof of the lemma).

For the first step, we define $g_0 = \abs{f}$; that is, for $u \in
B$, define $g_0(u) = \abs{f(u)}$. Note that $g_0$ is a
homomorphism from $B$ to $\Z$. Since $\min_{u \in \Gamma} f(u) =
0$, it follows that $g_0|_{\Gamma} = f|_{\Gamma}$. Furthermore,
every $u \in B$ admits $g_0(u) \geq 0$. Thus, $g_0$ has the two
properties described above.

At the $i\th$ step ($i >0$), define $g_i$ as follows
$$\forall \ u \in B \ \  g_i(u) = \left\{ \begin{array}{lr}
                    g_{i-1}(u) & u \not\in B_{r-i} \\
                    g_{i-1}(u) & u \in B_{r-i} \textrm{ and } g_{i-1}(u)
                    \neq i-1 \\
                    i+1 & u \in B_{r-i} \textrm{ and } g_{i-1}(u) =
                    i-1
                    \end{array} \right. $$
Since $\Gamma \cap B_{r-i} = \emptyset$, by induction
$g_i|_{\Gamma} = g_{i-1}|_{\Gamma} = f|_{\Gamma}$. Let $u$ be a
vertex in $B_{r-i}$ such that $g_{i-1}(u) = i-1$. Let $w$ be a
vertex in $B$ such that $w \sim u$. Then, $w \in B_{r-(i-1)}$.
Hence, by induction, $g_{i-1}(w) \geq i-1$. So, since $g_{i-1}$ is
a homomorphism, $g_{i-1}(w) = i$, which implies $g_i(u) - g_i(w) =
i+1 - g_{i-1}(w) = 1$. Thus, $g_i$ is indeed a homomorphism.
Furthermore, $\min_{u \in B_{r-i}} g_i(u) \geq i$.  So $g_i$
satisfies the properties described above.
\end{proof}

\subsection{Proof of Theorem \ref{lower bound theorem}} \label{sec: proof of lower bound gen}
Fix $n \in \N$. Set $G = G_n$, $f=f_n$, $v_0 = v_n$, $r = r(n)$,
and $S = V_n(r)$. We recall the following definitions. For $v \in
G$, we denote by $B_r(v)$ the ball of radius $r$ centered at $v$.
We say that \emph{$B_r(v)$ is of exact radius $r$}, if there
exists $u \in B_r(v)$ such that the distance between $v$ and $u$
in $G$ is at least $r$.  The following claim describes the size of
a collection of pairwise disjoint balls of exact radius  $r$ in
$G$.

\begin{clm} \label{clm: pairwise disjoint}
Let $V \subseteq G$ be a set of vertices of size $\abs{V} = k$.
Then, there exists a set $U \subseteq V$ of size $\abs{U} \geq
\lfloor {k/S^2} \rfloor$ such that
\begin{enumerate}
\item For all $u \in U$, the ball $B_r(u)$ is
of exact radius  $r$.
\item For all $u \neq u' \in U$, we have $B_r(u) \cap B_r(u') =
\emptyset$.
\end{enumerate}
\end{clm}

\begin{proof}
We prove the claim by induction on the size of $V$. Induction
base: If $\abs{V} < S^2$, then there is nothing to prove.
Induction step: Assume $\abs{V} \geq S^2$. If $r=0$, then set $U =
V$, and the claim follows. If $r>0$, then $S>1$, which implies
$\abs{V} > S$. Then, since $S$ is the size of the maximal ball of
radius $r$ in $G$, there exist $v,v' \in V$ such that the distance
between $v$ and $v'$ in $G$ is at least $r$. Thus, $B_r(v)$ is of
exact radius  $r$.

Denote $$B = \bigcup_{w \in B_r(v)} B_r(w).$$ Then, $\abs{B} \leq
S^2$. Denote $V' = V \setminus B$. Then, $\abs{V'} \geq k - S^2$
and $\abs{V'} < \abs{V}$. By induction, there exists a set $U'
\subseteq V'$ of size $\abs{U'} \geq \lfloor {\abs{V'}/S^2}
\rfloor \geq \lfloor{k /S^2} \rfloor - 1$ such that
\begin{enumerate}
\item For all $u \in U'$, the ball $B_r(u)$ is
of exact radius  $r$.
\item For all $u \neq u' \in U'$, we have $B_r(u) \cap B_r(u') =
\emptyset$.
\end{enumerate}
Set $U = U' \cup \set{v}$. So, $U \subseteq V$ of size $\abs{U}
\geq \lfloor {k/S^2} \rfloor$. To complete the proof it remains to
show that for all $u \in U'$, we have $B_r(v) \cap B_r(u) =
\emptyset$. Indeed, let $u \in U'$. Then, since $u \not \in B$,
the distance between $u$ and $B_r(v)$ in $G$ is more than $r$.
Thus, $B_r(v) \cap B_r(u) = \emptyset$.
\end{proof}

Returning to the proof of Theorem \ref{lower bound theorem}: By
Claim~\ref{clm: pairwise disjoint}, let $k = \lfloor
{\abs{G}}/{S^2} \rfloor-1$, and let $B_1,\ldots,B_k$ be a
collection of pairwise disjoint balls of exact radius  $r$ in $G$
such that for every $i \in [k]$, we have $v_0 \not\in B_i$. Note
that every $i \in [k]$ admits $\abs{B_i} \leq S$.

Let $i \in [k]$, and let $g \in \Hom_{v_0}(G \setminus B_i,\Z)$ be
a homomorphism that can be extended to a homomorphism in
$\Hom_{v_0}(G,\Z)$. Denote by $A_i$ the event $\set{|f(B_i)| \leq
r}$ and by $E_{i,g}$ the event $\set{f \big|_{G \setminus B_i} =
g}$. Since $g$ can be extended to a homomorphism in
$\Hom_{v_0}(G,\Z)$, we have $E_{i,g} \neq \emptyset$. Since there
are at most $2^{\abs{B_i}} \leq 2^{S}$ homomorphisms $f' \in
\Hom_{v_0}(G,\Z)$ such that $f'$ agrees with $g$ on $G \setminus
B_i$, we have
$$\abs{E_{i,g}} \leq 2^{S}.$$ Since $E_{i,g} \neq \emptyset$,
since $B_i$ is of exact radius  $r$, and since $v_0 \not \in B_i$,
using Lemma~\ref{boundary}, there exists a homomorphism $h \in
\Hom_{v_0}(G,\Z)$ such that $\abs{h(B_i)} \geq r+1$ and $h$ agrees
with $g$ on $G \setminus B_i$. That is, $h \not \in A_i$ and $h
\in E_{i,g}$, which implies $\abs{A_i \cap E_{i,g}} \leq
\abs{E_{i,g}}-1$. Hence,
\begin{eqnarray} \Pr\left[A_i \ \Big| \ E_{i,g} \right] =
\frac{\abs{A_i \cap E_{i,g}}}{\abs{E_{i,g}}} \leq
\frac{\abs{E_{i,g}}-1}{\abs{E_{i,g}}} \leq 1 - 2^{-S}. \label{eqn:
cc1}
\end{eqnarray} Note that for $j \neq i$, since $B_j \subseteq G \setminus
B_i$, the event $E_{i,g}$ determines $A_j$; that is, $E_{i,g} \cap
A_j$ is either $E_{i,g}$ or empty. Since (\ref{eqn: cc1}) holds
for any $g$ such that $E_{i,g} \neq \emptyset$,
\begin{eqnarray} \nonumber \Pr \bgive{ A_i}{ A_j : \ 1 \leq j < i } & = & \sum_g
\Pr\left[A_i \ \Big| \ A_j : \ 1 \leq j < i, \ E_{i,g} \right] \Pr
\left[E_{i,g} | A_j : \ 1 \leq j < i \right] \\ \nonumber & = &
\sum_g \Pr\left[A_i | E_{i,g} \right] \Pr \left[E_{i,g} | A_j : \
1 \leq j < i \right]
\\ & \leq & 1 - 2^{-S}, \label{eqn: cc2} \end{eqnarray}
where the sum is over all homomorphisms $g \in \Hom_{v_0}(G
\setminus B_i,\Z)$ such that $E_{i,g} \cap \set{A_j : 1 \leq j <
i}~\neq~\emptyset$ (which implies $E_{i,g} \cap \set{A_j : 1 \leq
j < i}=E_{i,g}$). Since (\ref{eqn: cc2}) holds for every $i \in
[k]$, and since $S \leq c \log \abs{G}$ (where $c < 1$),
\begin{eqnarray*}
    \Pr \Bracket{ R(f) \leq r} & \leq & \Pr [ \forall \ i \in [k] \ A_i ] \\
    & = & \prod_{i=1}^k \Pr \bgive{ A_i}{A_j : \ 1 \leq j < i } \\
    & \leq & \Soger{ 1 - 2^{-S} }^k  \leq \exp \Soger{ -
    \frac{k}{2^{S} } } \\ & \leq & e^2 \exp \Soger{- \frac{\abs{G}^{1-c}}{c^2 \log^2 \abs{G}}} =
o(1),
\end{eqnarray*}
where the last equality holds, since $\lim_{n \to \infty}
\abs{G_n} = \infty$. \qed

\section{The Cycle - $C_{n,k}$} \label{sec: the cycle}

In this section we study the range of a random homomorphism of the
graph $C_{n,k}$, where $n,k \in \N$, and $n$ is even. We consider
the graph $C_{n,k}$ mainly for two reasons: First, for logarithmic
$k$, the graph $C_{n,k}$ has almost linear diameter (the diameter
of $C_{n,k}$ is $\Omega(n)$, while the size of $C_{n,k}$ is $O(n
\log n)$), and still the range of a random homomorphism of
$C_{n,k}$ is constant (the range is $3$). Second, $C_{n,k}$ is
both vertex transitive and edge transitive.

The graph $C_{n,k}$ is a cycle of $n$ layers. Each layer has $k$
vertices, and is connected to both its adjacent layers by a
complete bi-partite graph. Thus, the degree of $C_{n,k}$ is $2k$.
Formally, the vertex set of $C_{n,k}$ is $\Z_n \times [k]$, and
the edge set of $C_{n,k}$ is defined by the relations $$\forall \
i \in \Z_n \ \ s,t \in [k] \ \ \ (i,s) \sim (i+1,t),$$ where
addition is modulo $n$. ($C_{n,k}$ is also the tensor product of
the $n$-cycle and the complete graph on $k$ vertices with
self-loops.) Denote by $\H_{n,k} = \Hom_{(0,1)}(C_{n,k},\Z)$, the
set of homomorphisms from $C_{n,k}$ to $\Z$ that map $(0,1)$ to
$0$. Since $n$ is even, $C_{n,k}$ is bi-partite, which implies
that $\H_{n,k} \neq \emptyset$.

We show a threshold phenomena (with respect to $k$) concerning the
range of a random homomorphism from $C_{n,k}$ to $\Z$. More
precisely, for $k(n) = 2 \log n + \omega(1)$, the range of a
random homomorphism is at most $3$, with high probability, and on
the other hand for $k(n) = 2 \log n - \omega(1)$, the range of a
random homomorphism is super-constant, with high probability. The
following two theorems make the above statements precise.

\begin{thm} \label{thm: upper bound for cnk}
Let $n \in \N$ be even, and let $k = k(n) = 2 \log n + \psi(n)$,
where $\psi: \N \to \mathbb{R}^+$ is such that $\lim_{n \to
\infty} \psi(n) = \infty$. Let $f_n \in_R \H_{n,k}$ be a random
homomorphism. Then
$$\Pr[R(f_n) \leq 3] = 1 - o(1). $$
\end{thm}

\begin{thm} \label{thm: lower bound for cnk}
Let $n \in \N$ be even, and let $k = k(n) = 2 \log n - \psi(n)$,
where $\psi: \N \to \mathbb{R}^+$ is monotone and $\lim_{n \to
\infty} \psi(n) = \infty$. Let $f_n \in_R \H_{n,k}$ be a random
homomorphism. Then
$$\Pr \Bracket{R(f_n) \geq \frac{2^{\psi(n-2)/4}}{\psi(n)}} = 1 - o(1).$$
\end{thm}

For the rest of this section we prove the above theorems. The
proof of Theorem~\ref{thm: upper bound for cnk} is deferred to
Section~\ref{sec: proof of upper bound}, and the proof of
Theorem~\ref{thm: lower bound for cnk} is deferred to
Section~\ref{sec: proof of lower bound}. We note that for $k=1$,
we have that $C_{n,k}$ is the $n$-cycle. Thus, $f_n \in_R
\H_{n,1}$ is a random walk bridge of length $n$ (a random walk
bridge is a random walk conditioned on returning to $0$). In this
case, Theorem~\ref{thm: lower bound for cnk} gives the bound
$$ \Pr \Bracket{ R(f_n) \geq \Omega \left( \frac{\sqrt{n}}{\log n} \right) } = 1-o(1). $$
This is consistent with the range of a random walk bridge (see
also Remark~\ref{rem: k is one}).

\subsection{Definitions} \label{sec: def of cnk}

Let $n,k \in \N$, where $n$ is even. For $i \in \Z_n$, the
\emph{$i$-layer} in $C_{n,k}$ is the set of vertices $\set{i}
\times [k]$. Recall that $\H_{n,k} = \Hom_{(0,1)}(C_{n,k},\Z)$ is
the set of homomorphisms from $C_{n,k}$ to $\Z$ that map $(0,1)$
to $0$. Denote by $\H^0_{n,k}$ the set of homomorphisms from
$C_{n,k}$ to $\Z$ that map the $0$-layer to $0$; that is,
$$\H_{n,k}^0 = \setb{f \in \H_{n,k}}{f(\set{0}\times[k]) = \set{0}}.$$
When $n$ and $k$ are clear we use $\H^0 = \H^0_{n,k}$.

For $f \in \H_{n,k}$ and $i \in \Z_n$, we say that the $i$-layer
is \emph{constant} in $f$, if $f$ gets the same value on the
entire $i$-layer; that is, $|f(\set{i}\times[k])|=1$. We say that
the $i$-layer is \emph{non-constant} in $f$, if $f$ gets different
values on the $i$-layer; that is, $|f(\set{i}\times[k])|>1$.
Define $\mathcal{NC}(f)$ to be the set of non-constant layers in
$f$; that is,
$$\mathcal{NC}(f) = \setb{i \in \Z_n}{|f(\set{i}\times[k])|>1}.$$
For $\ell \in [0,n]$, define $\H^0(\ell) = \H_{n,k}^0(\ell)$ to be
the set of homomorphisms in $\H^0$ that have exactly $\ell$
non-constant layers.

Loosely speaking, a homomorphism of $C_{n,k}$ corresponds to a
path on $\Z$ that starts at $0$ and ends at $0$. This motivates
the following definition. For an even integer $m \in \N$, denote
by $\P(m)$ the set of paths of length $m$ on $\Z$ that start at
$0$ and end at $0$; that is,
$$\P(m) = \setb{(S_0,S_1,\ldots,S_m) \in \Z^m}{\forall \ i \in [m] \ \ \abs{S_{i}-S_{i-1}} = 1
\textrm{ and } S_0 = S_m = 0}.$$ Note that $\abs{\P(m)} = {m
\choose m/2}.$

Consider the values of a homomorphism on the $1$-layer.  Since all
vertices in the $1$-layer are connected to a vertex that is mapped
to $0$, the value of a homomorphism on the $1$-layer corresponds
to a vector in $\set{1,-1}^k$.  In fact, it turns out that the
value of a homomorphism on a non-constant layer corresponds to a
$\set{1,-1}^k$ non-constant vector. Thus, the following definition
will be useful. Define
$$V = V_k = \set{1,-1}^k \setminus \set{(1,1,\ldots,1),(-1,-1,\ldots,-1)}.$$

\subsection{The Constant Layers}

In this section we show some properties of the constant layers.
First, we show that homomorphisms in $\H_{n,k}$ do not have two
adjacent non-constant layers. Second, we show that if the
$0$-layer is non-constant in a homomorphism $f \in \H_{n,k}$, then
we can think of $f$ as a homomorphism in $\H^0_{n-2,k}$. Third, we
show that, conditioned on a specific set of $\ell$ non-constant
layers, a random homomorphism in $\H^0_{n,k}$ corresponds to a
random walk bridge of length $n-2\ell$ (i.e., a random walk of
length $n-2\ell$ on $\Z$ that starts at $0$ and ends at $0$).

\subsubsection{No Two Adjacent Non-constant Layers}

\begin{clm} \label{clm: no two adjacent nonconst layer}
Let $f \in \H_{n,k}$ be a homomorphism. Assume that $i \in \Z_n$
is such that the $i$-layer is non-constant in $f$. Then there
exists $z \in \Z$ such that
$$f(\set{i+1} \times [k])= f(\set{i-1} \times [k]) = \set{z}.$$
In particular, both the $(i+1)$-layer and the $(i-1)$-layer are
constant in $f$.
\end{clm}

\begin{proof}
We prove the claim for the $(i+1)$-layer. The proof for the
$(i-1)$-layer is similar. Since the $i$-layer is non-constant in
$f$, there exist $s,t \in [k]$ such that $f(i,s) < f(i,t)$. Recall
that every $q \in [k]$ admits $(i+1,q) \sim (i,s)$ and $(i+1,q)
\sim (i,t)$. Thus, every $q \in [k]$ admits $f(i+1,q) =
f(i,s)+1=f(i,t)-1$. Setting $z = f(i,s)+1$, the claim follows.
\end{proof}

\subsubsection{If the $0$-layer is non-constant in $f$, we can think of $f$ as a homomorphism
of a smaller graph}

Let $f \in \H_{n,k} \setminus \H_{n,k}^0$ be a homomorphism. That
is, the $0$-layer is non-constant in $f$. Define $f_{\dwn}$ as
follows:$$ \forall \ i \in \Z_{n-2} \ \ s \in [k] \ \ \
f_{\dwn}(i,s) = f(i+1,s)- f(1,1) . $$

\begin{clm}
\label{clm: fdwn is homomo} Let $f \in \H_{n,k} \setminus
\H_{n,k}^0$ be a homomorphism. Then $f_\dwn$ is a homomorphism in
$\H^0_{n-2,k}$.
\end{clm}

\begin{proof}
Since the $0$-layer is non-constant in $f$, by Claim~\ref{clm: no
two adjacent nonconst layer}, there exists $z \in \set{1,-1}$ such
that \begin{eqnarray} f(\set{1} \times [k]) = f(\set{n-1} \times
[k]) = \set{z}. \label{eq: ab} \end{eqnarray} First, we show that
$f_\dwn$ is a homomorphism of $C_{n-2,k}$. Indeed, for all $i \in
[0,n-4]$ and for all $s,t \in [k]$, we have $f_\dwn(i+1,s) -
f_\dwn(i,t) = f(i+2,s)-f(i+1,t) \in \set{1,-1}$. Furthermore, for
all $s,t \in [k]$, by (\ref{eq: ab}), we have
$$f_\dwn(0,s) - f_\dwn(n-3,t) = f(1,s) - f(n-2,t) =  f(n-1,1)-f(n-2,t) \in \set{1,-1}.$$
Second, for all $s \in [k]$, we have $f_\dwn(0,s) = f(1,s) -
f(1,1) = z-z = 0$.
\end{proof}

In fact, the following claim holds.

\begin{clm} \label{clm: 0 non const f down}
Let $f \in_R \H_{n,k} \setminus \H_{n,k}^0$ be a random
homomorphism. Then $f_{\dwn}$ is uniformly distributed in
$\H_{n-2,k}^0$.
\end{clm}

\begin{proof}
We will show that the mapping
$$ \textrm{ from } \quad \H_{n,k} \setminus \H_{n,k}^0  \quad \textrm{ to } \quad \H_{n-2,k}^0 \times
\set{1,-1} \times \sr{ \set{0,2}^{k-1} \setminus \set{(0,\ldots,0)} } $$ %
defined by
$$ f \mapsto \sr{ f_\dwn,f(1,1),f(1,1) \cdot f(0,2),f(1,1)
\cdot f(0,3),\ldots,f(1,1) \cdot f(0,k) } $$ %
is a bijection, where
$f_\dwn \in \H_{n-2,k}^0$ is the homomorphism defined above,
$f(1,1) \in \set{1,-1}$, and $(f(1,1) \cdot f(0,2),f(1,1) \cdot
f(0,3),\ldots, f(1,1) \cdot f(0,k)) \in \set{0,2}^{k-1}$ is a
non-zero vector.

First, the mapping is injective. Indeed, let $f^1 \neq f^2$ be two
homomorphisms in $\H_{n,k} \setminus \H_{n,k}^0$. If $f^1(1,1)
\neq f^2(1,1)$, then the images of $f^1$ and $f^2$ are different
(in the second coordinate). Otherwise, assume that
\begin{eqnarray} f^1(1,1) = f^2(1,1). \label{eq: ab1}
\end{eqnarray}

There exist $i \in \Z_n$ and $s \in [k]$ such that $f^1(i,s) \neq
f^2(i,s)$. If either $i = 1$ or $i = n-1$, since the $0$-layer is
non-constant in both $f^1$ and $f^2$, using Claim~\ref{clm: no two
adjacent nonconst layer}, then $f^1(1,1) = f^1(i,s) \neq f^2(i,s)
= f^2(1,1)$ (contradicting (\ref{eq: ab1})). Otherwise, if $i =
0$, then $f^1(1,1) \cdot f^1(0,s) \neq f^2(1,1) \cdot f^2(0,s)$,
for $s >1$ (since $f^1(0,1) =f^2(0,1) = 0$), implying that the
images of $f^1$ and $f^2$ are different. Finally, if $i \in
[2,n-2]$, we have $f^1_\dwn(i-1,s) = f^1(i,s) - f^1(1,1) \neq
f^2(i,s) - f^2(1,1) = f^2_\dwn(i-1,s)$, so the images of $f^1$ and
$f^2$ are different (in the first coordinate).

Second, the mapping is surjective. Indeed, given a homomorphism $g
\in \H_{n-2,k}^0$, an integer $z \in \set{1,-1}$, and a non-zero
vector $(v_2,\ldots,v_k) \in \set{0,2}^{k-1}$, define
$$\forall \ i \in \Z_n \ s \in [k] \ \
f(i,s) = \left\{ \begin{array}{cc}
  0 & i = 0, s = 1 \\
  z \cdot v_s   & i = 0, s \neq 1 \\
  z & i = n-1 \\
  g(i-1,s)+z &  i \in [1,n-2].
\end{array} \right. $$
Thus, for every $i \in \Z_n$ and $s,t \in [k]$, (recall that
$g(\set{0} \times [k]) = \set{0}$),
$$f(i+1,s) - f(i,t) = \left\{\begin{array}{cc}
  g(0,s)+z-0 = z \in \set{1,-1}& i = 0,t=1\\
  g(0,s)+z-z \cdot v_t \in \set{1,-1}& i = 0,t\neq1\\
  0 - z \in \set{1,-1}& i = n-1, s = 1 \\
  z \cdot v_s - z \in \set{1,-1} & i  = n-1,s\neq1 \\
  z - (g(n-3,t)+z) \in \set{1,-1} & i = n-2 \\
  g(i,s)+z-(g(i-1,t)+z) \in \set{1,-1} & i \in [1,n-3],
\end{array} \right.$$
which implies $f \in \H_{n,k}$. Furthermore, since
$(v_2,\ldots,v_k)$ is a non-zero vector, it follows that
$f(\set{0} \times[k]) = \set{0,2z}$, which implies $f \not \in
\H^0_{n,k}.$ Finally, we will show that $f \mapsto
(g,z,v_2,\ldots,v_k).$ Indeed, for all $i \in \Z_{n-2}$ and $s \in
[k]$, we have $f_\dwn(i,s) = f(i+1,s) - f(1,1) =  g(i,s)+ z -
(g(0,1)+z) = g(i,s)$. Also $f(1,1) = g(0,1)+z = z$, and for all $s
\in [2,k]$, we have $f(1,1) \cdot f(0,s) = z \cdot z \cdot v_s =
v_s$.

The size of the range of the above defined mapping is
$\abs{\H_{n-2,k}^0} \cdot 2 \cdot (2^{k-1}-1)$. Therefore, for
every $g \in \H^0_{n-2,k}$,
$$ \Pr \Bracket{ f_{\dwn} = g } = \frac{ 2 \cdot (2^{k-1}-1) }{
\abs{ \H_{n,k} \setminus \H_{n,k}^0 } } =
\frac{1}{\abs{\H_{n-2,k}^0}}.
$$
\end{proof}

\subsubsection{Conditioned on the Set of Non-constant Layers, a Random Homomorphism is a
Random Walk Bridge}

Let $$I = \set{i_1 < \cdots < i_\ell} \subseteq [n-1]$$ be a set
of size $\ell$ such that for every $i \in [n-2]$, either $i \not
\in I$ or $i+1 \not \in I$ (or both). We think of $I$ as a set of
non-constant layers (recall Claim~\ref{clm: no two adjacent
nonconst layer}).


Denote by $\H(I,n)$ the set of homomorphisms $f$ in $\H^0_{n,k}$
such that $I$ is the set of non-constant layers in $f$ (we think
of $k$ as fixed). Recall that $\P(n-2\ell)$ is the set of paths on
$\Z$ of length $n-2\ell$ that start at $0$ and end at $0$, and
recall that
$$V = \set{1,-1}^k \setminus \set{(1,1,\ldots,1),(-1,-1,\ldots,-1)}.$$
For a homomorphism $f \in \H_{n,k}$, define the \emph{range of the
constant layers in $f$} to be
$$ \mathbf{RC}(f) = \setb{ f(i,1) }{ i \in \Z_n \textrm{ is such that the $i$-layer is constant in } f }.$$
For a path $(S_0,S_1,\ldots,S_{n-2\ell})$ in $P(n-2\ell)$, define
the \emph{range of the path} to be
$$ \Rng(S_0,S_1,\ldots,S_{n-2\ell}) = \setb{ S_i }{ 0 \leq i
\leq n-2\ell } . $$

Loosely speaking, the following proposition shows that,
conditioned on the set of non-constant layers, a random
homomorphism in $\H^0$ is a random walk bridge.

\begin{prop} \label{prop: given
nonconst the biject} Let $I = \set{ i_1 < i_2 < \cdots < i_\ell }
\subseteq [n-1]$ be a set of size $\ell$ such that for all $i \in
[n-2]$, either $i \not\in I$ or $i+1 \not\in I$. Then there exists
a bijection $\varphi$ between $\H(I,n)$ and $P(n-2\ell) \times
V^{\ell}$. Furthermore, denote $\varphi = (\varphi^1, \varphi^2)$.
Then for all $f \in \H(I,n)$,
$$ \mathbf{RC}(f) = \Rng(\varphi^1(f)) . $$
\end{prop}

For the rest of this section we prove Proposition~\ref{prop: given
nonconst the biject}.

\subsubsection{Proof of Proposition~\ref{prop: given nonconst the
biject}}

We prove the proposition by induction on $\ell$. The induction
step is based on the following claim. The claim shows that given a
non-constant layer in a homomorphism $f$ of $C_{n,k}$, we can
think of $f$ as a homomorphism of $C_{n-2,k}$.


In what follows, for simplicity, we use the following convention:
For a homomorphism $f \in \H_{n,k}$, and integers $i \in \N$ and
$s \in [k]$, we define $f(i,s) = f({i \pmod n},s)$.

\begin{clm} \label{clm: given nonconst the bij}
Let $I = \set{i_1 < \cdots < i_\ell} \subseteq [n-1]$ be a set of
size $\ell$ such that for every $i \in [n-2]$, either $i \not \in
I$ or $i+1 \not \in I$.  Let $f \in \H(I,n)$ be a homomorphism.
Define $f' \in \H(I \setminus \set{i_\ell}, n-2)$ by
$$ \forall \ i \in [0,n-3] \ s \in [k] \quad f'(i,s) = \left\{
\begin{array}{lr}
    f(i,s) & i < i_{\ell} \\
    f(i+2,s) & i \geq i_{\ell} .
\end{array} \right. $$
Define $v_f \in V$ by
$$ \forall \ s \in [k] \quad v_f(s) = f(i_{\ell},s) -
f(i_{\ell}-1,s) . $$ %
Then the map $\Phi:\H(I,n) \to \H(I \setminus \set{i_{\ell}}, n-2)
\times V$ defined by $\Phi(f) = (f',v_f)$ is a bijection.
\end{clm}

\begin{proof}
First we show that $f'$ and $v_f$ are well defined:

For $f'$, choose some $i \in [0,n-3]$ and $s,t \in [k]$. If $i+1 <
i_\ell$, then
$$ f'(i+1,s) - f'(i,t) = f(i+1,s) - f(i,t) \in \set{ 1,-1} . $$
If $i \geq i_{\ell}$ then
$$ f'(i+1,s) - f'(i,t) = f(i+3,s) - f(i+2,t) \in \set{ 1,-1} . $$
We are left with the case $i = i_\ell-1$.  Since the
$i_\ell$-layer is non-constant in $f$, by Claim~\ref{clm: no two
adjacent nonconst layer}, $f(i+2,t) = f(i,t)$.  Thus,
$$ f'(i+1,s) - f'(i,t) = f(i+3,s) - f(i,t) = f(i+3,s) - f(i+2,t)
\in \set{1,-1} . $$ %
So $f' \in \H_{n-2,k}$.  Since $i_{\ell} > 0$, for all $s \in
[k]$, we have $f'(0,s) = f(0,s) = 0$.  Thus, $f' \in
\H^0_{n-2,k}$. Also, for any layer $i < i_\ell$, since $f'(\set{i}
\times [k]) = f(\set{i} \times [k])$, we get that for any $i <
i_\ell$,
$$ \textrm{the $i$-layer is constant in } f' \quad \Leftrightarrow \quad
\textrm{the $i$-layer is constant in } f . $$ %
For any $i \geq i_{\ell}$, we have that $f'(\set{i} \times [k]) =
f(\set{i+2} \times [k])$.  Since $i+2 > i_\ell$, the $(i+2)$-layer
is constant in $f$, which implies that the $i$-layer is constant
in $f'$. So, the set of non-constant layers in $f'$ is the set $I
\setminus \set{i_\ell}$, and $f' \in \H(I \setminus
\set{i_{\ell}},n-2)$.

Now we show that $v_f$ is well defined:  Since the $i_\ell$-layer
is non-constant in $f$, using Claim~\ref{clm: no two adjacent
nonconst layer}, there exist $s,t \in [k]$ such that $v_f(s) \neq
v_f(t)$, so $v_f$ is in $V$.

To show that $\Phi$ is a bijection, we provide the inverse map
$\Psi = \Phi^{-1}$. Define $$\Psi : {\H(I \setminus
\set{i_\ell},n-2) \times V} \to \H(I,n)$$ as follows: For a pair
$g' \in \H(I \setminus \set{i_\ell},n-2)$ and $v \in V$, define $g
\in \H(I,n)$ by
$$ \forall \ i \in [0,n-1] \ s \in [k] \quad g(i,s) = \left\{
\begin{array}{lr}
    g'(i,s) & i < i_{\ell} \\
    g'(i-1,s) + v(s) & i = i_\ell \\
    g'(i-2,s) & i \geq i_{\ell} +1 ,
\end{array} \right. $$
and define $\Psi(g',v) = g$.

\Claim $g \in \H(I,n)$.

\begin{proof}
Choose some $i \in [0,n-1]$ and $s,t \in [k]$.  If $i+1 < i_\ell$,
then
$$ g(i+1,s) - g(i,t) = g'(i+1,s) - g'(i,t) \in \set{1,-1} . $$
If $i \geq i_\ell + 1$, then
$$ g(i+1,s) - g(i,t) = g'(i-1,s) - g'(i-2,t) \in \set{1,-1} . $$
If $i = i_\ell -1$, then, using Claim~\ref{clm: no two adjacent
nonconst layer},
$$ g(i+1,s) - g(i,t) = g'(i,s) + v(s) - g'(i,t) = v(s) \in
\set{1,-1} . $$ %
If $i = i_\ell$, then, using Claim~\ref{clm: no two adjacent
nonconst layer},
$$ g(i+1,s) - g(i,t) = g'(i-1,s) - (g'(i-1,t) + v(t)) \in \set{1,-1} . $$
So $g \in \H_{n,k}$.

Now, since $i_\ell > 0$, we have that $g(\set{0} \times [k]) =
g'(\set{0} \times [k]) = \set{0}$.  So $g \in \H_{n,k}^0$.

Finally, for $i < i_\ell$, we have that $g(\set{i} \times [k]) =
g'(\set{i} \times [k])$.  Thus, for $i < i_\ell$, we have that
$$ \textrm{the $i$-layer is constant in } g \quad \Leftrightarrow \quad
\textrm{the $i$-layer is constant in } g' . $$ %

If $i > i_\ell$, then $g(\set{i} \times [k]) = g'(\set{i-2} \times
[k])$. Since $i-2 > i_\ell - 2$, we have $i-2 \not \in I \setminus
\set{i_\ell}$. Hence, the $(i-2)$-layer is constant in $g'$, and
so the $i$-layer is constant in $g$. If $i = i_\ell$, then, since
$i-1 \not \in I \setminus \set{i_\ell}$, the $(i-1)$-layer is
constant in $g'$. Also
$$ v \not\in \set{ (1,1,\ldots,1) , (-1,-1,\ldots,-1) } . $$ %
So we have that the $i$-layer is non-constant in $g$ (since
$g(\set{i} \times [k]) = \set{g'(i-1,1)-1,g'(i-1,1)+1}$).

Thus, the set of non-constant layers of $g$ is
$$ (I \setminus \set{i_\ell}) \cup \set{ i_\ell}  = I
. $$ %
So $g \in \H(I,n)$, proving the claim.
\end{proof}

Now, let $f \in \H(I,n)$ and let $g = \Psi ( \Phi(f) )$.  Let $i
\in [0,n-1]$ and let $s \in [k]$. If $i-2 = i_\ell-1$, then since
the $i_\ell$-layer is non-constant in $f$, using Claim~\ref{clm:
no two adjacent nonconst layer},
$$ f'(i-2,s) = f(i_\ell-1,s) = f(i_\ell+1,s) = f(i,s) . $$
If $i-2 > i_\ell-1$, then $f'(i-2,s) = f(i,s)$.  So, if $i-2 \geq
i_\ell-1$, we have that $f'(i-2,s)= f(i,s)$. Thus, for all $i \in
[0,n-1]$ and $s \in [k]$,
$$  g(i,s) = \left\{
\begin{array}{lr}
    f'(i,s) = f(i,s) & i < i_{\ell} \\
    f'(i-1,s) + v_f(s) =
    f(i-1,s) + f(i_\ell,s) - f(i_\ell-1,s) & i = i_\ell \\
    f'(i-2,s) = f(i,s) & i \geq i_{\ell} +1 .
\end{array} \right. $$
So $g \equiv f$, and $\Psi = \Phi^{-1}$.  Thus, $\Phi$ is a
bijection as claimed.
\end{proof}

\begin{clm} \label{clm: RC eq RC}
Let $\Phi:\H(I,n) \to \H(I \setminus \set{i_\ell},n-2) \times V$
be the bijection defined in Claim~\ref{clm: given nonconst the
bij}. Let $f \in \H(I,n)$ be a homomorphism, and let $(f',v_f) =
\Phi(f)$. Then
$$ \mathbf{RC}(f) = \mathbf{RC}(f') . $$
\end{clm}

\begin{proof}
Let $z \in \mathbf{RC}(f)$.  So there exists $0 \leq i \leq n-1$
such that $f(i,1) = z$ and the $i$-layer is constant in $f$.

If $i < i_\ell$, then $f'(\set{i} \times [k]) = f(\set{i} \times
[k]) = \set{z}$, and so the $i$-layer is constant in $f'$.  Thus,
$z = f'(i,1) \in \mathbf{RC}(f')$.

We can exclude the case $i=i_\ell$, because the $i_\ell$-layer is
non-constant in $f$.

If $i = i_\ell+1$, then using Claim~\ref{clm: no two adjacent
nonconst layer}, for all $s \in [k]$, since the $i$-layer is
constant in $f$,
$$ f'(i-2,s) = f(i_\ell-1,s) = f(i_\ell+1,s) = f(i,s) = z .
$$ %
So, the $(i-2)$-layer is constant in $f'$, and $z = f'(i-2,1) \in
\mathbf{RC}(f')$.

If $i > i_\ell+1$, then $f'(\set{i-2} \times [k]) = f(\set{i}
\times [k]) = \set{z}$.  So, the $(i-2)$-layer is constant in
$f'$, and $z = f'(i-2,1) \in \mathbf{RC}(f')$.

This establishes $\mathbf{RC}(f) \subseteq \mathbf{RC}(f')$.

Now let $z \in \mathbf{RC}(f')$.  So, there exists $0 \leq i \leq
n-3$ such that $f'(i,1) = z$ and the $i$-layer is constant in
$f'$.

If $i < i_\ell$, then $f(\set{i} \times [k]) = f'(\set{i} \times
[k]) = \set{z}$.  So, the $i$-layer is constant in $f$, and $z =
f(i,1) \in \mathbf{RC}(f)$.

If $i \geq i_\ell$, then $f(\set{i+2} \times [k]) = f'(\set{i}
\times [k]) = \set{z}$.  So, the $(i+2)$-layer is constant in $f$
and $z = f(i+2,1) \in \mathbf{RC}(f)$.

So, $\mathbf{RC}(f') \subseteq \mathbf{RC}(f)$, which implies
$\mathbf{RC}(f) = \mathbf{RC}(f')$.
\end{proof}

Back to the proof of Proposition~\ref{prop: given nonconst the
biject}. Loosely speaking, we define $\varphi$ to be $\ell$
compositions of $\Phi$, where $\Phi$ is the map defined in
Claim~\ref{clm: given nonconst the bij}.

Let $f \in \H(I,n)$ be a homomorphism. For an integer $j \in
[1,\ell]$, define the set $$I_j = \set{i_1 < i_2 < \cdots <
i_{\ell-j}}.$$ Define $\ell$ homomorphisms $f_1,\ldots,f_\ell$ and
$\ell$ vectors $v_1,\ldots,v_\ell$ inductively as follows:
$(f_1,v_1) = \Phi(f)$, and $(f_j,v_j) = \Phi(f_{j-1})$. By
Claim~\ref{clm: given nonconst the bij}, for every $j \in
[1,\ell]$, we have $f_j \in \H(I_j,n-2j)$. Since $I_\ell =
\emptyset$, we have $f_\ell \in \H(\emptyset,n-2\ell)$. By
Claims~\ref{clm: given nonconst the bij} and ~\ref{clm: RC eq RC},
the map
$$f \mapsto (f_\ell,v_1,\ldots,v_\ell)$$ is a bijection, and
\begin{eqnarray} \mathbf{RC}(f)= \mathbf{RC}(f_1) = \cdots =
\mathbf{RC}(f_\ell). \label{eqn: ab4} \end{eqnarray}

We claim that there exists a bijection $\psi:
\H(\emptyset,n-2\ell) \to P(n-2\ell)$ such that for all $g \in
\H(\emptyset,n-2\ell)$, \begin{eqnarray} \mathbf{RC}(g) =
\Rng(\psi(g)). \label{eqn: ab3} \end{eqnarray} Indeed, let $g \in
\H(\emptyset,n-2\ell)$.  Then, $\psi(g) =
(g(0,1),g(1,1),\ldots,g(n-2\ell,1))$ is a path of length $n-2\ell$
in $\Z$, starting at $0$ and ending at $0$.  Furthermore,
(\ref{eqn: ab3}) holds. Note that homomorphisms in
$\H(\emptyset,n-2\ell)$ do not have non-constant layers.  So given
a path $(S_0,S_1,\ldots,S_{n-2\ell})$ in $P(n-2\ell)$, we can
define $g' \in \H(\emptyset,n-2\ell)$ by
$$ \forall \ i \in \Z_{n-2\ell} \ s \in [k] \quad g'(i,s) = S_i .
$$ %
Since $\psi(g') = (S_0,S_1,\ldots,S_{n-2\ell})$, it follows that
$\psi$ is a bijection.

Finally, define $\varphi$ as follows:
$$ \varphi:\H(I,n) \to P(n-2\ell) \times V^\ell \qquad
\varphi(f) = ((S_0,\ldots,S_{n-2\ell}),(v_1,\ldots,v_\ell)) , $$ %
where $(S_0,\ldots,S_{n-2\ell}) = \psi(f_\ell)$. By (\ref{eqn:
ab4}) and (\ref{eqn: ab3}), we have $\mathbf{RC}(f) =
\Rng(S_0,\ldots,S_{n-2\ell})$. \qed

\subsection{The Size of $\H^0$}

Fix $n,k \in \N$ ($n$ is even) for the rest of this section. We
consider $\H^0 = \H^0_{n,k}$. Note that, by Claim~\ref{clm: no two
adjacent nonconst layer}, every $f \in \H^0$ has at most $n/2$
non-constant layers. Thus,
\begin{eqnarray} \abs{\H^0} = \sum_{\ell=0}^{n/2}
\abs{\H^0(\ell)}, \label{eqn: h0 as sum} \end{eqnarray} where
$\H^0(\ell)$ is the set of homomorphisms of $C_{n,k}$ that have
exactly $\ell$ non-constant layers.

The following lemma gives a formula for the size of $\H^0(\ell)$.

\begin{lem} \label{lem: size of hnk0} Let $\ell \in
[0,n/2]$. Then
$$\abs{\H^0(\ell)} = {n  - \ell \choose \ell}
{n-2\ell \choose n/2 -\ell} (2^k-2)^\ell.$$
\end{lem}

\begin{proof} If $k = 1$, then
$$\H^0(\ell) = \left\{ \begin{array}{cc}
  {n \choose n/2} & \ell = 0 \\
  0 & \ell > 0 ,
\end{array} \right. $$ proving the lemma. So assume that $k>1$.
Define a family of sets $$\mathcal{I} = \setb{I \subseteq
[n-1]}{\abs{I} = \ell \textrm{ and for every $i \in [n-2]$, either
$i \not \in I$ or $i+1 \not \in I$}}.$$ Let $f \in \H^0(\ell)$ be
a homomorphism. Recall that $\mathcal{NC}(f)$ is the set of
non-constant layers of $f$. Since $f \in \H^0$, we have $0 \not
\in \mathcal{NC}(f)$, which implies that $\mathcal{NC}(f)
\subseteq [n-1]$. Since $f \in \H^0(\ell)$, we have
$\abs{\mathcal{NC}(f)} = \ell$. Using Claim~\ref{clm: no two
adjacent nonconst layer}, for every $i \in [n-2]$, either $i \not
\in \mathcal{NC}(f)$ or $i+1 \not \in \mathcal{NC}(f)$. Therefore,
$\mathcal{NC}(f) \in \mathcal{I}$. Furthermore, for every set $I
\in \mathcal{I}$, there exists a homomorphism $g \in \H^0(\ell)$
such that $\mathcal{NC}(g) = I$ (since $k>1$). Hence,
\begin{eqnarray} \mathcal I = \setb{\mathcal{NC}(f)}{f \in
\H^0(\ell)}. \label{eqn: aaa2}
\end{eqnarray} Define a map $\rho: \mathcal I \to {[n-\ell] \choose \ell}$
(where ${[n-\ell] \choose \ell}$ is the family of subsets of
$[n-\ell]$ of size $\ell$) as follows:
$$\forall \ I = \set{i_1 < \cdots< i_\ell} \in \mathcal{I} \ \ \rho(I) = \set{i_1 < i_2-1 < i_3-2
< \cdots < i_\ell - (\ell -1)}.$$ For $I = \set{i_1 < \cdots <
i_\ell} \in \mathcal{I}$, since $\abs{I} = \ell$ and since for
every $i \in [n-2]$, either $i \not \in I$ or $i+1 \not \in I$,
the set $\rho(I)$ is a subset of $[n-\ell]$ of size $\ell$. So
$\rho$ is well defined.
\begin{clm} \label{clm: I to I'}
The map $\rho$ is a bijection between $\mathcal{I}$ and ${[n-\ell]
\choose \ell}$.
\end{clm}
\begin{proof} For a set $J = \set{j_1 <
\cdots < j_\ell} \subseteq [n-\ell]$, define the map $\rho^{-1}$
by
$$\rho^{-1}(J) = \set{j_1 < j_2 + 1 < \cdots < j_\ell +
(\ell-1)}.$$ Thus, $\rho^{-1}(J)$ is of size $\ell$ and for every
$i \in [n-2]$, either $i \not \in \rho^{-1}(J)$ or $i+1 \not \in
\rho^{-1}(J)$, which implies $\rho^{-1}(J) \in \mathcal I$. Since
every $I \in \mathcal I$ admits $\rho^{-1}( \rho (I)) = I$, it
follows that $\rho$ is a bijection.
\end{proof}

Back to the proof of Lemma~\ref{lem: size of hnk0}. By (\ref{eqn:
aaa2}),
$$\abs{\H^0(\ell)} = \sum_{I \in \mathcal{I}} \abs{\H(I,n)},$$
where $\H(I,n)$ is the set of homomorphisms $f \in \H^0$ such that
$\mathcal{NC}(f) = I$. By Proposition~\ref{prop: given nonconst
the biject}, for every $I \in \mathcal{I}$,
$$\abs{\H(I,n)} = \abs{\P(n-2\ell)} \abs{V}^\ell = {n-2\ell
\choose n/2 - \ell} (2^k-2)^\ell.$$ By Claim~\ref{clm: I to I'},
$$\abs{\mathcal{I}} = {n-\ell \choose \ell}.$$ So the lemma follows.
\end{proof}

\subsection{Upper Bound for $k = 2 \log n + \omega(1)$
\label{sec: proof of upper bound} - Theorem~\ref{thm: upper bound
for cnk}}

In this part we show that for $k = 2 \log n + \omega(1)$, the
range of a random homomorphism from $C_{n,k}$ to $\Z$ is $3$, with
high probability. We use the formula for $\abs{\H^0(\ell)}$ to
conclude that $f$ has $n/2$ non-constant layers, with high
probability. Which implies that $f$ is ``almost'' constant.

\subsubsection{Many Non-constant Layers}

To prove Theorem~\ref{thm: upper bound for cnk}, we use the
following lemma, which states that there are $n/2$ non-constant
layers in a random homomorphism of $C_{n,k}$. Note that, by
Claim~\ref{clm: no two adjacent nonconst layer}, the maximal
number of non-constant layers in every homomorphism of $C_{n,k}$
is $n/2$.

\begin{lem} \label{lem: k geq 2logn many non-constant layer}
Let $n \in \N$ be even, and let $k = k(n) \geq 2 \log n +
\psi(n)$, where $\psi: \N \to \mathbb{R}^+$ is such that $\lim_{n
\to \infty} \psi(n) = \infty$. Let $f_n \in_R \H^0_{n,k}$ be a
random homomorphism. Then
$$\Pr[\abs{\mathcal{NC}(f_n)} = n/2] = 1-o(1).$$
\end{lem}

\begin{proof} Fix $n \in \N$, let $\H^0 = \H^0_{n,k(n)}$,
and let $f = f_n \in_R \H^0$ be a random homomorphism. For every
$\ell \in [0,n/2]$, denote $h^0(\ell) = \abs{\H^0(\ell)}$. By
Lemma~\ref{lem: size of hnk0}, every $\ell \in [0,n/2]$ admits
$$h^0(\ell) = {n  - \ell \choose \ell}
{n-2\ell \choose n/2-\ell} (2^k-2)^\ell,$$ which implies
\begin{eqnarray*} \frac{h^0(\ell+1)}{h^0(\ell)} =
\frac{(n-2\ell)(n-2\ell-1)(n/2-\ell)^2(2^k-2)}{(n-\ell)(\ell+1)(n-2\ell)(n-2\ell-1)}
= \frac{(n/2-\ell)^2(2^k-2)}{(n-\ell)(\ell+1)}. \end{eqnarray*}
Thus, since $k \geq 2\log n + \psi(n)$ and since $n \geq 2$, every
$0 \leq \ell < n/2$ admits
$$\frac{h^0(\ell+1)}{h^0(\ell)} \geq 2^{\psi(n)}.$$ Thus, every $0 \leq
\ell < n/2$ admits $$\frac{h^0(n/2)}{h^0(\ell)} \geq
(2^{\psi(n)})^{n/2 - \ell},$$ which implies
$$\sum_{0 \leq \ell < n/2 } h^0(\ell) =
\sum_{0 \leq \ell < n/2 } h^0(n/2) \frac{h^0(\ell)}{h^0(n/2)} \leq
h^0(n/2) \sum_{0 \leq \ell < n/2} (2^{\psi(n)})^{\ell-n/2} =
o(h^0(n/2)),$$ where the last equality holds, since $\lim_{n \to
\infty} \psi(n) = \infty$. Thus,
$$\Pr[\abs{\mathcal{NC}(f)} = n/2] =
\frac{h^0(n/2)}{\sum_{0 \leq \ell \leq n/2 } h^0(\ell)} = 1 -
o(1).$$
\end{proof}

\subsubsection{Proof of Theorem~\ref{thm: upper bound for cnk}}

Denote $f = f_n$, and consider the following two cases:

{\bf Case one:} Assume that $f \in_R \H^0_{n,k}$. By
Lemma~\ref{lem: k geq 2logn many non-constant layer}, with
probability $1 - o(1)$, there are $n/2$ non-constant layers in
$f$. Thus, since $n$ is even and since the $0$-layer is constant
in $f$, using Claim~\ref{clm: no two adjacent nonconst layer}, all
the odd layers are non-constant in $f$. Therefore, all the even
layers are mapped to $0$. Hence, with probability $1 - o(1)$, we
have $f(C_{n,k}) \subseteq \set{-1,0,1}$, and $R(f) \leq 3$.

{\bf Case two:} Assume that $f \in_R \H_{n,k} \setminus
\H^0_{n,k}$. By Claim~\ref{clm: 0 non const f down}, we have that
$f_{\dwn}$ is uniformly distributed in $\H^0_{n-2,k}$. Hence, by
Lemma~\ref{lem: k geq 2logn many non-constant layer}, with
probability $1 - o(1)$, there are $n/2 -1$ non-constant layers in
$f_{\dwn}$. Thus, by definition of $f_{\dwn}$, including the
$0$-layer, there are $n/2$ non-constant layers in $f$. Hence, by
Claim~\ref{clm: no two adjacent nonconst layer}, all the even
layers are non-constant in $f$, and all the odd layers are
constant in $f$. Hence, with probability $1 - o(1)$, either
$f(C_{n,k}) \subseteq \set{0,1,2}$ or $f(C_{n,k}) \subseteq
\set{0,-1,-2}$. Therefore, with probability $1 - o(1)$, we have
$R(f) \leq 3$. \qed

\subsection{The Number of Non-constant Layers Determines the Range}

In the previous section we have seen that if the number of
non-constant layers is large, then the range is small. In this
section we show that, in fact, the number of non-constant layers
determines the range of a random homomorphism. The following lemma
gives a lower bound on the range of a random homomorphism of
$C_{n,k}$ with exactly $\ell$ non-constant layers. The lower bound
is determined by $\ell$. We note that a similar upper bound can be
proven.

\begin{lem} \label{lem: unfix
det range} Let $n \in \N$ be even, and let $k = k(n) \in \N$. Let
$\ell = \ell(n) \in \N$ be such that $\lim_{n \to \infty}
(n-2\ell) = \infty$. Let $f_n \in_R \H^0_{n,k}(\ell)$ be a random
homomorphism from $C_{n,k}$ to $\Z$ with exactly $\ell$
non-constant layers such that $f(\set{0}\times[k])=\set{0}$. Then
for every $\alpha > 0$,
$$\Pr \left[R(f_n) \geq \alpha \sqrt{n-2\ell}  \right] \geq (1-o(1))(1-2\alpha^2),$$
where the $o(1)$ term is as $n$ tends to infinity, and is
independent of $\alpha$.
\end{lem}

Loosely speaking, the proof of the lemma is as follows.
Conditioned on $f$ having $\ell$ non-constant layers, $f$
corresponds to a random walk bridge of length $n-2\ell$. Since the
range of such a walk is roughly $\sqrt{n-2\ell}$, the range of $f$
is roughly $\sqrt{n-2\ell}$.

\subsubsection{Proof of Lemma~\ref{lem: unfix det range}}

Recall that for an even $m \in \N$, we defined $\P(m)$ to be the
set of paths on $\Z$ of length $m$ that start at $0$ and end at
$0$. The following proposition shows that, with high probability,
the range of a random walk bridge of length $m$ is at least
$\Omega(\sqrt{m})$.

\begin{prop} \label{prop: range of bridge}
Let $m \in \N$ be even, and let $(S_0,S_1,\ldots,S_{m}) \in_R
\P(m)$ be a random walk bridge of length $m$. Then for every
$\alpha > 0$,
$$\Pr[ \abs{\set{S_0,S_1,\ldots,S_{m}}} \geq \alpha \sqrt{m}] \geq
(1 - o(1)) ( 1-2\alpha^2 ),$$ where the $o(1)$ term is as $m$
tends to infinity, and is independent of $\alpha$.
\end{prop}

First we use the proposition to prove the lemma. Fix $n \in \N$.
Partition $\H^0(\ell)$ as follows:
$$\H^0(\ell) = \bigcup_{I} \H(I,n),$$
where $I \subseteq [n-1]$, and $\H(I,n)$ is the set of
homomorphisms $f \in \H^0(\ell)$ such that $\mathcal{NC}(f) = I$.
Note that by Claim~\ref{clm: no two adjacent nonconst layer}, if
$\H(I,n) \neq \emptyset$, then $\abs{I} = \ell$ and for all $i \in
[n-2]$, either $i \not\in I$ or $i+1 \not\in I$.

Fix $I$ such that $\H(I,n) \neq \emptyset$, and let $g \in_R
\H(I,n)$. Denote
$$((S_0,S_1,\ldots,S_{n-2\ell}),(v_1,\ldots,v_\ell)) =
\varphi(g),$$ where $\varphi$ is the bijection given by
Proposition~\ref{prop: given nonconst the biject}. Thus, by
Proposition~\ref{prop: given nonconst the biject}, we have $R(g)
\geq \abs{\set{S_0,S_1,\ldots,S_{n-2\ell}}}$, and
$(S_0,S_1,\ldots,S_{n-2\ell})$ is uniformly distributed in
$\P(n-2\ell)$.

Let $f = f_n \in_R \H^0(\ell)$ be a random homomorphism. Then, for
all $I$ such that $\H(I,n) \neq \emptyset$,
$$ \Pr \left[ R(f) \geq \alpha \sqrt{n-2\ell} \ \Big|
f \in \H(I,n) \right] \geq \Pr \left[
\abs{\set{S_0,S_1,\ldots,S_{n-2\ell}}} \geq \alpha \sqrt{n-2\ell}
\right],$$ where  $(S_0,,\ldots,S_{n-2\ell}) \in_R \P(n-2\ell)$ is
a random walk bridge. Thus, we have
\begin{eqnarray} \nonumber \Pr \left[ R(f) \geq \alpha \sqrt{n-2\ell} \right]  &
= & \sum_{I}  \Pr \left[ R(f) \geq \alpha \sqrt{n-2\ell} \ \Big| f
\in \H(I,n) \right] \Pr[f \in \H(I,n)] \\
& \geq & \Pr \left[ \abs{\set{S_0,S_1,\ldots,S_{n-2\ell}}} \geq
\alpha \sqrt{n-2\ell} \right] \label{eqn: b1} ,
\end{eqnarray}
where the sum is over all sets $I \subseteq [n-1]$ such that
$\H(I,n) \neq \emptyset$.
Thus, by Proposition~\ref{prop: range of bridge}, since $\lim_{n
\to \infty} (n-2\ell) = \infty$,
$$(\ref{eqn: b1}) \geq (1-o(1)) \left(1-2\alpha^2 \right) , $$
and the $o(1)$ term is as $n$ tends to infinity, and is
independent of $\alpha$. \qed

\subsubsection{Proof of Proposition~\ref{prop: range of bridge}}

If $\alpha \geq 1$, then the Proposition follows. Thus, assume
$\alpha < 1$.

Let $T \in [m]$. Before proving the proposition we show that a
path in $\Z$ of length $m$ from $0$ to $0$ that passes through $T$
corresponds to a path in $\Z$ of length $m$ from $0$ to $2T$.
Formally,

\begin{clm} \label{clm: random bridge to T is as walk to 2T}
There exists a bijection between paths in $\P(m)$ that pass
through $T$ and paths in $\Z$ of length $m$ that start at $0$ and
end at $2T$.
\end{clm}

\begin{proof} Let $(S_0,S_1,\ldots,S_{m}) \in \P(m)$ be such that
there exists $j \in [m]$ that admits $S_j = T$. Let $j^* = \min
\setb{j \in [m]}{S_j = T}.$ The bijection is reflecting the path
around $T$ from $j^*$ onwards. That is, for $j \in [0,j^*]$ set
$S'_j = S_j$, and for $j \in [j^*+1,m]$ set $S'_j = 2T - S_j$.
Thus, $S'_0 = 0$, $S'_{j^*} = T$ and $S'_{m} = 2T$. Furthermore,
$(S'_0,S'_1,\ldots,S'_{m})$ is a path in $\Z$ of length $m$ such
that $S'_0 = 0$ and $S'_{m} = 2T$. Note that $j^*$ is also the
first time that $S'$ passes through $T$.

Since every path in $\Z$ of length $m$ that starts at $0$ and ends
at $2T$ passes through $T$, the above defined map is a bijection.
Indeed, we show how to invert the above defined map. Let
$0=S'_0,\ldots,S'_{m} = 2T$ be a path in $\Z$ of length $m$. Let
$j^* = \min \setb{j \in [m]}{S'_j = T}$. For $j \in [0,j^*]$ set
$S_j = S'_j$, and for $j \in [j^*+1,m]$ set $S_j = 2T - S'_j$.
\end{proof}

Since there are ${m \choose m/2-T}$ paths in $\Z$ of length $m$
that start at $0$ and end at $2T$, using Claim~\ref{clm: random
bridge to T is as walk to 2T}, \begin{eqnarray} \Pr[
\abs{\set{S_0,S_1,\ldots,S_{m}}} \geq T ] \geq \Pr[ \exists \ j
\in [m] \ : \ S_j = T] = \frac{{m \choose m/2-T}}{{m \choose
m/2}}. \label{eqn: a1} \end{eqnarray} Using Stirling's formula
(recall that for $x \geq 0,$ we have $1-x \leq e^{-x}$ and $1+x
\leq e^x$), substituting $T = \alpha \sqrt{m}$, we have
\begin{eqnarray*} \nonumber (\ref{eqn: a1}) & = & \frac{(m/2)! (m/2)!}{(m/2-T)!(m/2+T)!} \\ \nonumber &
= & (1 - o(1)) \left(1- \frac{2T}{m} \right)^{-m/2+T} \left(1+
\frac{2T}{m} \right)^{-m/2-T} \\ \nonumber & = & (1 - o(1))
\left(1- \frac{4T^2}{m^2} \right)^{-m/2+T} \left(1+ \frac{2T}{m}
\right)^{-2T}
\\
\nonumber & \geq & (1 - o(1)) e^{\frac{4T^2}{m^2}(m/2-T) - \frac{4T^2}{m}} \\
&
\geq & (1 - o(1)) \left( 1- \frac{2T^2}{m} \right) \\
\nonumber & = & (1 - o(1)) ( 1-2\alpha^2 ),
\end{eqnarray*} where the $o(1)$ term
is as $m$ tends to infinity, and is independent of $\alpha$, since
$\alpha < 1$. \qed

\subsection{A Lower Bound for $k = 2 \log n - \omega(1)$ - Theorem~\ref{thm: lower bound for cnk}}
\label{sec: proof of lower bound}

In this part we show that for $k = 2 \log n - \omega(1)$, the
range of a random homomorphism from $C_{n,k}$ to $\Z$ is
super-constant, with high probability. The proof plan is as
follows: First, we prove that there are many constant layers in a
random homomorphism $f$. Second, using Lemma~\ref{lem: unfix det
range}, we conclude that the range of $f$ is large.

\subsubsection{Many Constant Layers}

The following lemma shows that a random homomorphism of $C_{n,k}$
has many constant layers.

\begin{lem} \label{lem: k leq 2logn many constant layer}
Let $n \in \N$ be even, and let $k = k(n) = 2 \log n - \psi(n)$,
where $\psi: \N \to \mathbb{R}^+$ is monotone and $\lim_{n \to
\infty} \psi(n) = \infty$. Let $f_n \in_R \H_{n,k}$ be a random
homomorphism. Let $\beta: \N \to \mathbb{R}^+$ be such that for
large enough $n \in \N$, we have $\beta(n) \leq n/4$. Then for
large enough $n \in \N$, we have
$$ \Pr \Bracket{ \abs{ {\cal NC}(f_n) } > n/2 - \beta } \leq 16
\beta^2 2^{-\psi(n)} . $$
\end{lem}

\begin{proof} Fix some large enough $n \in \N$ such that $\beta \leq n/4$.
For $\ell \in [0,n/2]$, set $h^0(\ell) = \abs{\H_{n,k}^0(\ell)}$.
Let $\ell \in [n/2 - \beta,n/2-1]$, then by Lemma~\ref{lem: size
of hnk0}, since $\beta \leq n/4$,
$$ \frac{ h^0(\ell+1) }{ h^0(\ell) } = \frac{(n/2 - \ell)^2
(2^k-2) }{ (n-\ell) (\ell+1) } \leq \frac{ 4 \beta^2 2^k }{ n (n -
2\beta) } \leq 8 \beta^2 2^{-\psi(n)} . $$ %
Hence, setting $\gamma = 8 \beta^2 2^{-\psi(n)}$, every $\ell \in
[n/2 - \beta,n/2]$ admits
$$ \frac{h^0(\ell) }{ h^0(n/2-\beta) } \leq \gamma^{\ell -
(n/2-\beta)} . $$ %
Thus, (if $\gamma \geq 1/2$, the claim follows) since $\gamma <
1/2$,
\begin{eqnarray*}
    \Pr \Bracket{ \abs{ {\cal NC}(f_n) } > n/2 - \beta } & \leq &
    \frac{1}{ h^0(n/2-\beta) }  \sum_{\ell = n/2 - \beta +
    1}^{n/2} h^0(\ell)
    \leq  \sum_{i=1}^{\infty} \gamma^i \leq 2 \gamma .
\end{eqnarray*}
\end{proof}

\subsubsection{Proof of Theorem~\ref{thm: lower bound for cnk}}

First we consider homomorphisms in $\H^0$.

\begin{clm}\label{clm: f in H0 lower}
Let $n \in \N$ be even, and let $k = k(n) = 2 \log n - \psi(n)$,
where $\psi: \N \to \mathbb{R}^+$ is monotone and $\lim_{n \to
\infty} \psi(n) = \infty$. Let $f_n \in_R \H^0_{n,k}$ be a random
homomorphism. Let $\eps: \N \to \mathbb{R}^+$ be such that for
large enough $n \in \N$, we have $\eps(n) \leq 1/8$ and $\lim_{n
\to \infty} \eps(n) 2^{\psi(n)/2} = \infty$. Then
$$ \Pr \Bracket{ R(f_n) \geq \sqrt{2} \ \eps  2^{\psi(n)/4} } \geq
(1-o(1)) (1-2\eps)^2 . $$ %
\end{clm}

\begin{proof} Consider large enough $n \in \N$ such that $\eps(n) \leq 1/8$.
Set $\beta = \beta(n) = \eps(n) 2^{\psi(n)/2}$. Note
that $\lim_{n \to \infty} \beta(n) = \infty$. Since $1 \leq k = 2
\log n - \psi(n)$ and $\eps < 1/4$, we have $\beta \leq n/4$. Let
$\ell = \ell(n) \in [0,n/2 - \beta]$, then $\lim_{n \to \infty}
(n-2\ell) \geq \lim_{n \to \infty} 2 \beta(n) = \infty$.  Thus, by
Lemma~\ref{lem: unfix det range}, for any $\alpha
> 0$,
$$ \Pr \bgive{ R(f_n) \geq \alpha \sqrt{2 \beta} }{ f_n \in
\H^0(\ell) } \geq \Pr \bgive{ R(f_n) \geq \alpha \sqrt{n-2\ell} }{
f_n \in \H^0(\ell) } \geq (1-o(1)) (1-2\alpha^2) . $$ %
Thus, for any $\alpha > 0$,
\begin{eqnarray} \label{eqn:  Pr [ R(f) geq alpha sqrt(2 beta) ]}
    \Pr \Bracket{ R(f_n) \geq \alpha \sqrt{2 \beta} } & \geq &
    \sum_{\ell = 0}^{n/2 - \beta} \Pr \bgive{ R(f_n) \geq \alpha
    \sqrt{2 \beta} }{ f_n \in \H^0(\ell) } \Pr \Bracket{
     f_n \in \H^0(\ell) } \nonumber \\
    & \geq & (1-o(1)) (1-2\alpha^2) \Pr \Bracket{ \abs{ {\cal
    NC}(f_n) } \leq n/2 - \beta } .
\end{eqnarray}
By Lemma~\ref{lem: k leq 2logn many constant layer},
$$ (\ref{eqn:  Pr [ R(f) geq alpha sqrt(2 beta) ]}) \geq (1-o(1))(1-2\alpha^2) (1-16
\beta^2 2^{-\psi(n)}) = (1-o(1)) (1-2\alpha^2) (1-16 \eps^2) . $$
Taking $\alpha = \sqrt{\eps}$, since $\eps \leq 1/8$, we have
$$ \Pr \Bracket{ R(f_n) \geq \sqrt{2} \ \eps 2^{\psi(n)/4} } \geq
(1-o(1)) (1-2\eps)(1-16 \eps^2) \geq (1-o(1)) (1-2\eps)^2 . $$ %
\end{proof}

\begin{rem} \label{rem: k is one} If $k=1$, then $f_n \in_R \H_{n,k}^0$ is
a random walk bridge of length $n$. In this case, Claim~\ref{clm:
f in H0 lower} gives the bound
$$ \Pr \Bracket{ R(f_n) \geq 2^{1/4} \eps \sqrt{n} } \geq (1-o(1))
(1-2\eps)^2 . $$ %
\end{rem}

Back to the proof of Theorem~\ref{thm: lower bound for cnk}. Set
$\eps = 2^{-1/2} \psi(n)^{-1}$. Note that $\lim_{n \to \infty}
\eps(n) = 0$, and $\lim_{n \to \infty} \eps(n) 2^{\psi(n)/2} =
\infty$. We denote $f = f_n \in_R \H_{n,k}$. Consider the
following two cases:

{\bf Case one:} Assume that $f \in_R \H^0_{n,k}$. By
Claim~\ref{clm: f in H0 lower}, with probability at least
$$(1-o(1)) (1-2\eps)^2 = 1-o(1),$$ the range of $f$ is at least
$$\psi(n)^{-1}2^{\psi(n)/4}.$$

{\bf Case two:} Assume that $f \in_R \H_{n,k} \setminus
\H^0_{n,k}$. By Claim~\ref{clm: 0 non const f down}, we have that
$f_{\dwn}$ is uniformly distributed in $\H^0_{n-2,k}$. Thus, by
Claim~\ref{clm: f in H0 lower}, with probability at least
$$(1-o(1)) (1-2\eps)^2 = 1-o(1),$$ the range of $f_{\dwn}$ is at least
$$\psi(n-2)^{-1}2^{\psi(n-2)/4}.$$ By definition of
$f_{\dwn}$, the size of the range of $f$ is at least the size of
the range of $f_{\dwn}$. So, with probability $1-o(1)$, the range
of $f$ is at least $\psi(n)^{-1}2^{\psi(n-2)/4}$ (recall that
$\psi(n)$ is monotone). \qed

\section{Further Research} \label{sec: further research}

We list some possible further research directions regarding random
homomorphisms of graphs:

\begin{enumerate}

\item Let $G$ be a graph. A function $f:G \to \Z$ is called
\emph{Lipschitz} if it satisfies
$$ \forall \ v \sim u \in G \quad \abs{ f(u) - f(v)} \leq 1 . $$
Note that a homomorphism is always Lipschitz, but typically the
set of homomorphisms is much smaller than the set of Lipschitz
functions.  For a graph $G$ and a vertex $v \in G$, define
$\Lip_v^0(G,\Z)$ to be the set of all Lipschitz functions from $G$
to $\Z$ that map $v$ to $0$.
\begin{conj*}
Let $\set{G_n}$ be a family of bi-partite graphs with $\lim_{n \to
\infty} \abs{G_n} = \infty$.  Assume that for all $n$, $G_n$ has
maximal degree $d$ ($d$ independent of $n$). Let $f_n \in_R
\Hom_{v_n}^0(G_n,\Z)$ be a random homomorphism, and let $g_n \in_R
\Lip_{v_n}^0(G_n,\Z)$ be a random Lipschitz function.  Then,
$$ \frac{\E \br{ R(f_n) }}{ \E \br{ R(g_n)} } = \Theta (1) ,$$
where $\Theta(\cdot)$ depends on $d$.
\end{conj*}

\item Let $G$ be a bi-partite graph, and let $\Delta$ be the diameter of $G$.
For any homomorphism $f \in \Hom_v^0(G,\Z)$, we have that $R(f) =
O(\Delta)$.  But this naive bound should not be the typical bound,
at least not for symmetric graphs.
\begin{conj*}
Let $\set{G_n}$ be a family of vertex transitive bi-partite graphs
with $\lim_{n \to \infty} \abs{G_n} = \infty$.  Let $\Delta_n$ be
the diameter of $G_n$, and let $f_n \in_R \Hom_{v_n}^0(G_n,\Z)$ be
a random homomorphism. Then,
$$ \E \br{ R(f_n) } = O \sr{ \sqrt{ \Delta_n } } . $$
\end{conj*}
Note that the conjecture is false if the assumption of vertex
transitivity is dropped.  Consider, for example, the star graph
with $\Delta$ arms of length $\Delta$.  This graph has expected
range of $\Theta(\Delta)$.

\item For $d \in \N$, let $\O(d)$ be the set of all
\emph{odd} positive integers at most $d$. For an even integer $n
\in \N$, let $G_{n,d}$ be the graph whose vertices are $\Z_n$, and
edges are defined by the relations
$$ i \sim j \quad \Longleftrightarrow \quad \exists \ d' \in \O(d) \ \ i=j + d' {\pmod
n} . $$ %
Note that $G_{n,d}$ is vertex transitive and bi-partite (for all
$n$ and $d$). Note also that the diameter of $G_{n,d}$ is
$\Theta(n/d)$, and the degree of $G_{n,d}$ is $\Theta(d)$.

\begin{conj*}
There exists a constant $c > 0$ such that
\begin{enumerate}
\item If $d = c \log(n) - \omega(1)$,
$$ \Pr \br{ R(f_n) \geq \omega(1) } = 1-o(1),$$
where $f_n \in_R \Hom_0^0(G_{n,d},\Z)$ is a random homomorphism.
\item There exists a constant $b \in \N$ such that if $d = c\log(n) + \omega(1)$, then
$$ \Pr \br{ R(f_n) \leq b } = 1-o(1),$$
where $f_n \in_R \Hom_0^0(G_{n,d},\Z)$ is a random homomorphism.
\end{enumerate}
\end{conj*}

\item In this paper we only consider homomorphisms of graphs into
$\Z$.  Instead of $\Z$ consider the infinite star with $k$ arms.
That is, the graph $S_k$, whose vertices are
$$ V(S_k) = \Set{ (i,s) }{ 0 < i \in \N \ , \ 1 \leq s \leq k } \cup \set{0 } , $$
and edges are $(1,s) \sim 0$ for all $1 \leq s \leq k$, and
$$ \forall \ 1 < i \in \N \ , \ 1 \leq s \leq k \ , \quad  (i,s) \sim (i \pm 1,
s) . $$ %
It can be shown \cite{Ori} that a random homomorphism from the
interval $[n]$ to $S_3$, has range $O(\log(n))$ (as opposed to
$\Theta(\sqrt{n})$ when the homomorphism is into $\Z = S_2$).  Is
this the case for all bi-partite graphs $G$?  That is, let
$\set{G_n}$ be a family of bi-partite graphs such that $\lim_{n
\to \infty} \abs{G_n} = \infty$. Let $f_n \in_R
\Hom_{v_n}^0(G_n,S_3)$.  Is it true that
$$ \E \br{ R(f_n) } = O \sr{ \log \abs{G_n} } ? $$
Or maybe even
$$ \lim_{n \to \infty} \Pr \br{ R(f_n) = O \sr{ \log \abs{G_n} }
} = 1 ? $$

\end{enumerate}

\end{document}